\documentclass[titlepage,twoside,12pt]{article}
\usepackage{amssymb}
\usepackage{amsfonts}
\begin{document}


\pagenumbering{arabic}
\setcounter{page}{1}

\pagenumbering{arabic}

{\LARGE \bf What Scalars Should We Use ? } \\ \\

{\bf Elem\'{e}r E Rosinger \\ Department of Mathematics \\ University of Pretoria \\
Pretoria, 0002 South Africa \\
e-mail : eerosinger@hotmail.com} \\ \\

{\bf Abstract} \\

There are compelling historical and mathematical reasons why we ended up, among others in
Physics, with using the scalars given by the real numbers in $\mathbb{R}$, or the complex
numbers in $\mathbb{C}$. Recently, however, infinitely many easy to construct and use {\it
other} algebras of scalars have quite naturally emerged in a number of branches of Applied
Mathematics. These algebras of scalars can deal with the long disturbing difficulties
encountered in Physics, related to such phenomena as "infinities in Physics",
"re-normalization", the "Feynman path integral", and so on. Specifically, as soon as one is
dealing with scalars in algebras which - unlike $\mathbb{R}$ and $\mathbb{C}$ - are {\it no
longer} Archimedean, one can deal with a large variety of "infinite" quantities and do so
within the usual rules and with the usual operations of algebra. Here we present typical
constructions of these recently emerged algebras of scalars, most of them
non-Archimedean. \\ \\

{\bf Algebraic Note} \\

In order to avoid possible undesirable overlap between terminology used in mathematics, and on
the other hand, in physics, we recall here the following. \\
In mathematics the notion of {\it field} means a {\it ring} in which every non-zero element
has an inverse with respect to multiplication. Typical {\it fields} are given by the set
$\mathbb{Q}$ of rational numbers, the set $\mathbb{R}$ of real numbers, the set $\mathbb{C}$
of complex numbers, as well as the set $\mathbb{H}$ of Hamilton's quaternions. \\
The more general notion of {\it ring} means an algebraic structure with addition "+", and with
multiplication ".", with the addition "+" having the properties of a commutative group, while
the multiplication "." being distributive with respect to the addition "+". We should note
that such a multiplication "." need not always be commutative. \\
An essential {\it difference} between fields, and on the other hand, rings in general is that
in a field a product is zero, if and only if at least one of the factors is zero. In rings it
can happen that products are zero, without any of the factors being zero. This happens, for
instance, with square matrices of order at least 2. \\
The notion of {\it algebra} over a given field $\mathbb{K}$ means a {\it ring} $A$ whose
elements can be multiplied by elements in the field $\mathbb{K}$, namely, if $\lambda \in
\mathbb{K}$ and $x \in A$, then $\lambda x \in A$. This multiplication also has certain
properties with respect to the addition and multiplication in the ring $A$. Typically, square
matrices of real numbers, and of order at least 2, are algebras over $\mathbb{R}$ which are
not commutative, neither are they fields, since not all such non-zero matrices have a
multiplicative inverse. Similar properties hold for square matrices of complex numbers. \\
By {\it scalars} we shall understand in the sequel either real or complex {\it numbers} in
$\mathbb{R}$, respectively in $\mathbb{C}$, or elements in various {\it more general
algebras}, elements of which are supposed to be used instead of usual numbers. \\ \\ \\

\begin{quote}
\begin{quote}
\begin{quote}

"... Whatever is thought by us is either conceived through itself, or involves the concept of
another ... \\
So one must either proceed to infinity, or all thoughts are resolved into those which are
conceived through themselves ... \\
Every idea is analyzed perfectly only when it is demonstrated a priori that it is
possible ... \\
Since, however, it is not in our power to demonstrate the possibility of things in a perfectly
a priori way, that is, to analyze them into God and nothing, it will be sufficient for us to
reduce their immense multitude to a few, whose possibility can either be supposed and
postulated, or proved by experience ..." \\

\hspace*{3cm} G W Leibniz \\
\hspace*{3cm} Of an Organum or \\
\hspace*{3cm} Ars Magna (circa 1679)

\end{quote}
\end{quote}
\end{quote}

\bigskip
{\bf 0. Infinite regression of thought is not possible} \\

Neither in Physics, nor in Mathematics, or for that matter, not in any other human made theory
is it possible to keep endlessly considering each concept involved as being defined in terms
or more basic ones. \\
Consequently, at any given moment in time, each specific theory must start with a number of
concepts which are taken for granted, based on a variety of reasons mostly outside of the
respective theory. \\
Much of present day Physics takes for granted in this way concepts such as {\it real number},
{\it complex number}, or for that matter, the {\it Differential and Integral Calculus} based
on them, as initially constructed by Newton and Leibniz in the late 1600s. \\

Recently, there have been a number of instances in studies of Physics where such traditional
basic concepts have been reconsidered and replaced with other ones. \\
For instance, in connection with Brownian motion, the Loeb measure and integration, based on
nonstandard real numbers, has proved to be particularly useful, see Cutland or Albeverio
et.al. See also De Leo \& Rotelli regarding quaternionic electro-weak theory. \\
Earlier studies, see Lambek [1,2], have shown that one may gain a useful insight by going
beyond the familiar real and complex numbers, and employing Hamilton's quaternions, introduced
back in the mid 1800s. \\
Yet more impressive have been recent studies in Quantum Gravity which show the usefulness of
stepping out of the traditional Cantorian set theory into the far more general theory of
Categories, see Isham. \\

In this way, what had for so long proved to be a number of {\it de facto unstated axioms}
concerning the traditional use of certain mathematical concepts or structures in Physics have
lately been reconsidered. Here however one should note that the main interests of physicists
are naturally not related to the issue of basic mathematical concepts. Consequently,
physicists tend early in their career to become accustomed to the use of certain such basic
mathematical concepts, and are not particularly prone to reconsider them time and again. \\
In this regard, the recent survey paper of Baez [1,2] is a welcome exception. It shows to what
extent it may be useful in several theories of Physics to go beyond the traditionally well
entrenched use of real and complex numbers as the only {\it scalars} ever to be
employed. \\ \\

{\bf 1. Scalars : Should we hold unaware to an Unstated Axiom ?} \\

The following citation usually attributed to R P Feynman :

\begin{quote}

"The whole purpose of Physics is to find a number, with decimal point, etc. ! \\
Otherwise you haven't done anything."

\end{quote}

may indeed express a well founded practical point of view. However, it in fact leaves open the
question :

\begin{quote}

What kind of number ?  \\
Or more generally, what kind of scalars ?

\end{quote}

And as things go, this question is never really asked explicitly, since the answer appears to
be there instantly and a priori, namely : a real or complex number, of course ... \\

Thus, in Quantum Mechanics, among other theories of Physics, we operate with the : \\

{\bf Unstated Axiom about Scalars}, or in short {\bf UAS} : \\

All numbers must belong to the field $\mathbb{R}$ of real scalars, or at most to the larger
field $\mathbb{C}$ of complex scalars. \\

\hfill $\Box$ \\

Certainly, the acceptance - knowingly, or otherwise - of such or any other axiom which
delimits, and in fact, restricts the kind of scalars to be used from then on exclusively,
can have considerable implications upon any theory of physics subsequently developed. Such an
acceptance in its consequences is indeed much unlike the fact whether the respective
physicist thinks, speaks and writes in English, Russian or Chinese, for instance. After all,
all civilized languages are more or less equivalent with one another when it comes to the
realms of meanings they can express. Not to mention that, unlike with poetry for instance, the
translation of scientific texts from one to another among such languages is quite accurate and
without particular difficulties. \\

Here, in order to help keeping up an appropriate wider perspective, it is worth recalling the
use in Quantum Field Theory of the rather strange {\it Grassmann anticommuting scalars}, of
the respective unusual integration based on them, and of the corresponding Grassmann path
integrals, see Zee [pp. 124-126]. \\
The presence in Quantum Theory of such "detuors" from the usual real or complex numbers, thus
from UAS, may indeed help maintaining a more appropriate perspective with respect to the
variety of possible scalars to be used. \\

Although this may at first sound rather off the mark, we shall show that :

\begin{itemize}

\item one can {\it easily construct} a large variety of algebras of scalars which extend the
usual field of complex numbers,

\end{itemize}

Needless to say, the introduction of such extended algebras of scalars may have further
implication in Quantum Theory, or in Physics in general. Indeed, two important features of
these algebras are that :

\begin{itemize}

\item they are not Archimedean,

\item they have zero divisors, thus they are not division algebras, in particular, they are
not fields, and cannot be included in fields, although they include the usual field of complex
numbers $\mathbb{C}$.

\end{itemize}

The fact that these extended algebras of scalars are not Archimedean may prove to represent a
{\it significant advantage} due to the fact that one does no longer hit in them so easily and
so often upon the so called "infinities in physics". Indeed, these algebras can contain
scalars which are "infinitely large" and which, nevertheless, can be subjected without any
exception to all the usual rules and operations in an algebra. \\

Here it should be noted that there has for more than four decades by now been a precedent in
this regard, namely, with the particular case of nonstandard reals $^*\mathbb{R}$. Indeed ,
the nonstandard reals are {\it not} Archimedean either, and they allow an easy treatment of a
rich variety of "infinitely large" elements. \\
However, the algebras of scalars we deal with in this paper are still larger, and also more
general then the field of nonstandard reals. Furthermore, in view of the presence in them of
zero divisors, these algebras are no longer fields, unlike the case with the nonstandard
reals. \\
More importantly yet, the algebras of scalars dealt with in this paper - although larger and
more general - can nevertheless be defined, studied and used through the usual, classical,
standard methods of algebra, unlike the nonstandard numbers which need a rather elaborate
theory that is at present not widely enough familiar. \\

As for the presence of zero divisors, the eventual occurrence of possible disadvantages coming
from such a situation is not so clear at present, and consequently, it may require further
studies. Indeed, as long as these extended algebras of scalars contain the usual field of
complex numbers and do not modify in any way the operations on such usual numbers, we may only
encounter the issue of zero divisors when we go beyond complex numbers. \\
It is however important to note here that we have by now for long been working, and at some
ease, with algebras which have zero divisor, namely, {\it algebras of usual matrices} of real
or complex numbers. \\

Lastly, we can note that the use of the algebras of scalars introduced in this paper need not
prevent one from the use of Calculus. Indeed, as shown recently, Differential and Integral
Calculus can naturally be extended to a large variety of such algebras, see Bertram
et.al. \\ \\

{\bf 2. Three Classical Fields of Scalars} \\

One of the consequences in Quantum Mechanics of the unstated axiom UAS mentioned in section 1,
is that, in Quantum Field theory, for instance, we soon hit so called "infinities" which lead
to all sorts of technical and not particularly convincing juggling, such as for instance
"renormalization", and which try to avoid them. \\
However, regardless of such or other possible difficulties, it may - after several millennia
of mathematics, and several centuries of scientific physics - occur to us that the above axiom
could, and in fact, should be subjected to a certain conscious and closer examination. \\

Let us try and present here certain arguments which may, or may not support that axiom. \\

We can start by noting that neither the field of real scalars in $\mathbb{R}$, nor of the
complex ones $\mathbb{C}$ happened to come our way by some chance. Indeed, these two fields
have properties which make them {\it unique}, and through those properties, also make them
quite natural to come to our attention. \\
For instance, it is well known since the 19th century that $\mathbb{R}$ is the {\it only}
field which is totally ordered, Archimedean, and also complete in its order. \\
In the 1930s, Pontrjagin proved that the only fields which are not totally disconnected
topologically are the real numbers $\mathbb{R}$, the complex numbers $\mathbb{C}$, and the
Hamiltonian quaternions $\mathbb{H}$. Furthermore, we have the strict inclusions \\

(2.1) $~~~~~~ \mathbb{R} ~\subset~ \mathbb{C} ~\subset~ \mathbb{H} $ \\

and each smaller field is a subfield of the larger one in which it happens to be included.
This result arose as an answer to a question asked by Kolmogorov and Pontrjagin, see also
Rosinger [1]. \\

With respect to multiplication, each of these three fields are {\it associative}, the first
two are also {\it commutative}, while $\mathbb{H}$ is noncommutative. On the other hand,
$\mathbb{R}$ is {\it not} algebraically closed, while the other two fields are. \\
In this way, the quaternions $\mathbb{H}$ may appear to take a second place when compared to
$\mathbb{R}$ and $\mathbb{C}$, owing to their noncommutativity, although there have been
studies showing certain significant advantages in case we would use them in the formulation of
equations of physics, see Horn for more recent literature, as well as the series of papers on
quaternions at www.arXiv.org \\

There is a further {\it uniqueness} property which is of interest here. Namely, $\mathbb{R},~
\mathbb{C},~ \mathbb{H}$ and the octonions $\mathbb{O}$ are the only {\it normed division
algebras}, Baez [1, p. 150]. \\
The concept of division algebra is explained in the sequel, while a normed algebra is any
algebra $A$ with a usual norm $||~||$, which in addition satisfies the condition $|| x~y || =
|| x ||~ || y ||$, for every $x, y \in A$. \\

The fact that the octonions, discovered in the mid 1800s not much after the quaternions
$\mathbb{H}$, did not elicit much interest in physics is most likely due to the fact that, as
an algebra, they fail not only to be commutative, but also associative, and this latter
failure makes computations with octonions rather cumbersome. Yet recently, there has been a
growing interest in the use of octonions in physics, namely, in Quantum Logic, Special
Relativity and Supersymmetry, Baez [1,2]. \\

Two features of the above fields of scalars $\mathbb{R},~ \mathbb{C}$ and $\mathbb{H}$ are
important to note. \\
As {\it vector spaces} on $\mathbb{R}$, they are ordered, and as such, they are Archimedean. \\
As for their multiplications, they are {\it without zero divisors}, or in other words, they
are {\it division algebras}, Baez [1,2]. \\

Let us look into these two properties in some more detail. \\

The order relation $\leq$ on $\mathbb{R}$ is the well known obvious one, and as such, it is a
{\it total} order, namely, for every $x, y \in \mathbb{R}$, we have $ x \leq y$, or $y \leq
x$. \\
The order relations $\leq$ on $\mathbb{C}$ and $\mathbb{H}$, respectively, are induced by the
order relation on $\mathbb{R}$, noting that as {\it vector spaces} on $\mathbb{R}$, we have
the isomorphisms \\

(2.2) $~~~~~~ \mathbb{C} ~\cong~ \mathbb{R}^2,~~~~ \mathbb{H} ~\cong~ \mathbb{R}^4 $ \\

and on every Euclidean space $\mathbb{R}^d$, with $d \geq 1$, we have the {\it partial order}
$\leq$ induced naturally by the total order $\leq$ on $\mathbb{R}$, according to \\

(2.3) $~~~~~~ ( x_1,~.~.~.~, x_d ) ~\leq~ ( y_1,~.~.~.~, y_d ) ~~~\Longleftrightarrow~~~
                    x_1 ~\leq~ y_1,~~.~~.~~.~~, x_d ~\leq~ y_d $ \\

About these partial orders, the important point to note is that each of them is {\it
Archimedean}, this property being shared by all orders in (2.3) on Euclidean spaces. Namely \\

(2.4) $~~~~~~ \begin{array}{l}
                     \exists~~~ u \in \mathbb{R}^d,~~ u ~\geq~ 0 ~: \\ \\
                     \forall~~~ x \in \mathbb{R}^d,~~ x ~\geq~ 0 ~: \\ \\
                     \exists~~~ n \in \mathbb{N} ~: \\ \\
                     ~~~~~ n~ u ~\geq~ x
               \end{array} $ \\

since we can obviously take $u = ( 1,~.~.~.~, 1 ) \in \mathbb{R}^d$. \\

The inevitable, and not seldom troubling consequence of this Archimedean property is that in
every Euclidean space, therefore, in view of (2.2), also in $\mathbb{R},~ \mathbb{C}$ and
$\mathbb{H}$, we shall have for the above $u$ \\

(2.5) $~~~~~~ u + u + u + ~.~.~.~ ~=~ \infty \notin \mathbb{R}^d $ \\

This is, therefore, the source of all those so called "infinities in physics", namely, the
tacit acceptance of the above axiom {\bf UAS} which, either we like it, or not, does in fact
{\it confine} us to Archimedean scalars. \\

By the way, within the field $^*\mathbb{R}$ of nonstandard real numbers, the relation (2.5)
does {\it not} hold, since this field is {\it not} Archimedean. Consequently, we have the
presence of {\it infinitely large} elements in $^*\mathbb{R}$, and we have for instance \\

(2.6) $~~~~~~ 1 + 1 + 1 + ~.~.~.~ ~=~ \omega \in~ ^*\mathbb{R} $ \\

with $\omega$ being a well defined infinitely large scalar, with which one can do without any
restriction whatsoever all the usual operations in a field. \\

When it comes to the second common property of $\mathbb{R},~ \mathbb{C}$ and $\mathbb{H}$,
namely, that each of them is without zero divisors, or equivalently, it is a division algebra,
this means that they satisfy the property \\

(2.7) $~~~~~~ \begin{array}{l}
                       \forall~~ x,~ y \in \mathbb{R} ~: \\ \\
                       ~~~~ x~. ~y ~=~ 0 ~~~\Longrightarrow~~~ x ~=~ 0 ~~\mbox{or}~~ y ~=~ 0
               \end{array} $ \\

and similarly for $\mathbb{C}$ and $\mathbb{H}$. \\

This property may often be useful, in particular, in simplifying relations involving products.
However, we are quite accustomed to work in algebras in which the multiplication does not have
this property, like for instance is the case with the multiplication of matrices of real or
complex numbers. \\ \\

{\bf 3. What other Scalars are there available ?} \\

It was observed back in the 1800s by Cayley and Dickson that the way one obtains $\mathbb{C}$
from $\mathbb{R}$, namely, by identifying each complex number $z \in \mathbb{C}$ with a pair
of real numbers $( a, b ) \in \mathbb{R} \times \mathbb{R}$, can be repeated in order to
obtain $\mathbb{H}$ from $\mathbb{C}$, and then $\mathbb{O}$ from $\mathbb{H}$. Furthermore,
one can repeat that construction indefinitely, Baez [1,2]. Unfortunately however, the scalars
one obtains are more and more losing their nicer properties, for instance, they fail even to
be associative, Baez [1,2]. \\
Also in the 1800s, Clifford constructed a variety of associative algebras which in the simpler
cases recall the Cayley-Dickson construction. \\
Recently in Davenport a four dimensional commutative algebra $\mathbb{D}$ was suggested
instead of the quaternions $\mathbb{H}$. \\
And if we are to keep more focused on physics, we can recall here again the use in Quantum
Field Theory of the rather strange Grassmann anticommuting scalars, mentioned earlier. \\

What is, however, somewhat overlooked when searching for algebras of scalars is the simple
fact that, starting with any field, or even merely an algebra $\mathbb{K}$, one can easily
construct {\it infinitely many} other algebras. Here, for the convenience of familiarity, we
shall assume that the algebras considered are over the field $\mathbb{R}$ of real numbers, or
the field $\mathbb{C}$ of complex numbers. However, the method of construction presented next
is automatically valid if one starts with every algebra $\mathbb{K}$, be it commutative or not,
associative or not. \\

This method of construction is well known in the branch of Mathematical Logic called Model
Theory, where in its more general form is one of the fundamental constructions going under the
name of {\it reduced products}, see Rothmaler. A particular case of such reduced products are
the {\it reduced powers} which we shall use in the sequel in the construction of large classes
of algebras of scalars, see (3.4). \\

By the way, a further particular case of reduced powers are the ultra powers used in one of
the customary constructions of the field $^*\mathbb{R}$ of nonstandard reals, see Albeverio
et.al., or Cutland. However, this ultra power construction proves to be somewhat too
particular for the purposes pursued in this paper. Not to mention that it involves one in the
technical complications of "transfer", "internal" or "external" sets, and so on, typical for
dealing with nonstandard numbers. \\

On the other hand, the use of reduced products and reduced powers can be done quite easily,
and within the framework of usual mathematics, more precisely, algebra, thus {\it without} the
need for any nonstandard methods. \\

Here it can be noted that reduced powers were extensively used in the construction of large
classes of algebras of generalized functions, as well as generalized scalars, both employed in
the solution of a large variety of nonlinear PDEs, see section 4, below. \\

The mentioned large class of algebras of scalars, constructed as reduced powers, is obtained
as follows. \\
Let $\mathbb{K}$ be any algebra, and $\Lambda$ be any infinite set, then we consider \\

(3.1) $~~~~~~ \mathbb{K}^\Lambda $ \\

which is the set of all functions $x : \Lambda \longrightarrow \mathbb{K}$, and which is
obviously an algebra with the point-wise operations on such functions. Furthermore, if
$\mathbb{K}$ is commutative or associative, the same will hold for the algebra
$\mathbb{K}^\Lambda$. \\

A possible problem with these algebras $\mathbb{K}^\Lambda$ is that they have {\it zero
divisors}, thus they are {\it not} division algebras, in the case of nontrivial $\mathbb{K}$,
that is, when $0,~ 1 \in \mathbb{K}$, and $1 \neq 0$. For instance, if we take $\mathbb{K} =
\mathbb{R}$ and $\Lambda = \mathbb{N}$, then \\

(3.2) $~~~~~~ \begin{array}{l}
       ( 1, 0, 0, 0,~.~.~.~ ),~ ( 0, 1, 0, 0,~.~.~.~ ) \in \mathbb{K}^\Lambda ~=~
                                       \mathbb{R}^\mathbb{N}, \\ \\
       ( 1, 0, 0, 0,~.~.~.~ ),~ ( 0, 1, 0, 0,~.~.~.~ ) ~\neq~
                                   0 \in \mathbb{K}^\Lambda ~=~ \mathbb{R}^\mathbb{N}
               \end{array} $ \\

and yet \\

(3.3) $~~~~~~ \begin{array}{l}
                       ( 1, 0, 0, 0,~.~.~.~ )~.~ ( 0, 1, 0, 0,~.~.~.~ ) ~=~
                             ( 0, 0, 0,~.~.~.~ ) ~=~ \\ \\
                        ~~~~=~ 0 \in \mathbb{K}^\Lambda ~=~ \mathbb{R}^\mathbb{N}
              \end{array} $ \\

This issue, however, can be dealt with in the following manner, to the extent that it may
prove not to be convenient. Let us take any {\it ideal} ${\cal I}$ in $\mathbb{K}^\Lambda$ and
construct the {\it quotient algebra} \\

(3.4) $~~~~~~ A ~=~ \mathbb{K}^\Lambda / {\cal I} $ \\

This quotient algebra construction has {\it four} useful features, namely

\begin{enumerate}

\item It allows for the construction of a large variety of algebras $A$ in (3.4),

\item The algebras $A$ in (3.4) are fields, if and only if the corresponding ideals ${\cal I}$
are {\it maximal} in $\mathbb{K}^\Lambda$,

\item The algebras $A$ in (3.4) are without zero divisors, thus are division algebras, if and
only if the corresponding ideals ${\cal I}$ are {\it prime} in $\mathbb{K}^\Lambda$,

\item If $\mathbb{K} = \mathbb{R}$, then there is a simple way to construct such ideals
${\cal I}$ in $\mathbb{K}^\Lambda$, due to their one-to-one correspondence with {\it filters}
on the respective infinite sets $\Lambda$.

\end{enumerate}

Here it should be mentioned that the above construction of quotient algebras in (3.4) has
among others a well known particular case, namely, the construction of the field
$^*\mathbb{R}$ of nonstandard reals, where one takes $\mathbb{K} = \mathbb{R}$, while
${\cal I}$ is a maximal ideal which corresponds to an ultrafilter on $\Lambda$. \\

Facts 2 and 3 above are a well known matter of undergraduate algebra. \\
Fact 4 is recalled here briefly for conveninece. First we recall that each $x \in
\mathbb{R}^\Lambda$ is actually a function $x : \Lambda \longrightarrow \mathbb{R}$. Let us
now associate with each $x \in \mathbb{R}^\Lambda$ its {\it zero set} given by $Z ( x ) ~=~
\{~ \lambda \in \Lambda ~|~ x ( \lambda ) ~=~ 0 ~\}$, which therefore is a subset of
$\Lambda$. \\

Further, let us recall the concept of {\it filter} on the set $\Lambda$. A family ${\cal F}$
of subsets of $\Lambda$, that is, a subset ${\cal F} \subseteq {\cal P} ( \Lambda )$, is
called a {\it filter} on $\Lambda$, if and only if it satisfies the following three
conditions \\

\bigskip
(3.5) $~~~~~~ \begin{array}{l} 1.~~ \phi \notin {\cal F} \neq \phi \\ \\
                               2.~~ J,~ K \in {\cal F} ~~\Longrightarrow~~ J \bigcap K \in
                               {\cal F} \\ \\
                               3.~~ \Lambda \supseteq K \supseteq J \in {\cal F}
                                              ~~\Longrightarrow~~ K \in {\cal F}
                         \end{array} $ \\

Given now an ideal ${\cal I}$ in $\mathbb{R}^\Lambda$, let us associate with it the set of
zero sets of its elements, namely \\

(3.6) $~~~~~~ {\cal F}_{\cal I} ~=~ \{~ Z ( x ) ~~|~~ x \in {\cal I} ~\}
                                                   ~\subseteq~ {\cal P} ( \Lambda ) $ \\

Then \\

(3.7) $~~~~~~ {\cal F}_{\cal I} ~~\mbox{is a filter on}~ \Lambda $ \\

\medskip
Indeed, ${\cal I} \neq \phi$, thus ${\cal F}_{\cal I} \neq \phi$. Further, assume that
$Z ( x ) = \phi$, for a certain $x \in {\cal I}$. Then $x ( \lambda ) \neq 0$, for $\lambda
\in \Lambda$. Therefore we can define $y : \Lambda ~\longrightarrow~ \mathbb{R}$, by
$y ( \lambda ) = 1 / x ( \lambda )$, with $\lambda \in \Lambda$. Then however $y . x = 1 \in
\mathbb{R}^\Lambda$, hence ${\cal I}$ cannot be an ideal in $\mathbb{R}^\Lambda$, which
contradicts the hypothesis. In this way condition 1 in (3.5) holds for ${\cal F}_{\cal I}$. \\
Let now $x,~ y \in {\cal I}$, then clearly $x^2 + y^2 \in {\cal I}$, and $Z ( x^2 + y^2 ) =
Z ( x ) \bigcap Z ( y )$, thus condition 2 in  (3.5) is also satisfied by
${\cal F}_{\cal I}$. \\
Finally, let $x \in {\cal I}$ and $K \subseteq \Lambda$, such that $K \supseteq Z ( x )$. Let
$y$ be the characteristic function of $\Lambda \setminus K$. Then $x . y \in {\cal I}$, since
${\cal I}$ is an ideal. But now obviously $Z ( x . y ) = K$, which shows that
${\cal F}_{\cal I}$ satisfies as well condition 3 in (3.5). \\

There is also the {\it converse} construction. Namely, let ${\cal F}$ be any filter on
$\Lambda$, and let us associate with it the set of functions \\

(3.8) $~~~~~~ {\cal I}_{\cal F} ~=~ \{~ x : \Lambda ~\longrightarrow~ \mathbb{R} ~~|~~ Z ( x )
                                         \in {\cal F} ~\} ~\subseteq~ \mathbb{R}^\Lambda $ \\

Then \\

(3.9) $~~~~~~ {\cal I}_{\cal F} ~~\mbox{is an ideal in}~ \mathbb{R}^\Lambda $ \\

Indeed, for $x, y \in \mathbb{R}^\Lambda$, we have $Z ( x + y ) \supseteq Z ( x ) \bigcap
Z ( y )$, thus $x,~ y \in {\cal I}_{\cal F}$ implies that $x + y \in {\cal I}_{\cal F}$. Also
$Z ( x . y ) \supseteq Z ( x )$, therefore $x \in {\cal I}_{\cal F},~ y \in
\mathbb{R}^\Lambda$ implies that $x . y \in {\cal I}_{\cal F}$. Further we note that
$Z ( c x ) = Z ( x )$, for $c \in \mathbb{R},~ c \neq 0$. Finally, it is clear that
${\cal I}_{\cal F} \neq \mathbb{R}^\Lambda$, since $x \in {\cal I}_{\cal F} ~\Longrightarrow~
Z ( x ) \neq \phi$, as ${\cal F}$ satisfies condition 1 in (3.5).Therefore (3.9) does indeed
hold. \\

Let now ${\cal I},~ {\cal J}$ be two ideals in $\mathbb{R}^\Lambda$, while
${\cal F},~ {\cal G}$ are two filters on $\Lambda$. Then it is easy to see that \\

(3.10) $~~~~~~ \begin{array}{l}
                       {\cal I} ~\subseteq~ {\cal J} ~~~\Longrightarrow~~ {\cal F}_{\cal I}
                               ~\subseteq~ {\cal F}_{\cal J} \\ \\
                       {\cal F} ~\subseteq~ {\cal G} ~~~\Longrightarrow~~ {\cal I}_{\cal F}
                       ~\subseteq~ {\cal I}_{\cal G}
                \end{array} $ \\

We can also note that, given an ideal ${\cal I}$ in $\mathbb{R}^\Lambda$ and a filter
${\cal F}$ on $\Lambda$, we have by iterating the above constructions in (3.6) and (3.8) \\

(3.11) \quad $ \begin{array}{l}
                 {\cal I} ~~~\longrightarrow~~~ {\cal F}_{\cal I} ~~~\longrightarrow~~~
                             {\cal I}_{{\cal F}_{\cal I}} ~=~ {\cal I} \\ \\
                 {\cal F} ~~~\longrightarrow~~~ {\cal I}_{\cal F} ~~~\longrightarrow~~~
                             {\cal F}_{{\cal I}_{\cal F}} ~=~ {\cal F}
                        \end{array} $ \\

Indeed, in view of (3.6), (3.8), we have for $s \in \mathbb{R}^\Lambda$ the equivalent
conditions \\

$~~~~~~ x \in {\cal I} ~~~\Longleftrightarrow~~ Z ( x ) \in {\cal F}_{\cal I}
                  ~~~\Longleftrightarrow~~~ x \in {\cal I}_{{\cal F}_{\cal I}} $ \\

Further, for $J \subseteq \Lambda$, we have the equivalent conditions \\

$~~~~~~ J \in {\cal F}_{{\cal I}_{\cal F}} ~~~\Longleftrightarrow~~~ J ~=~ Z ( s ),
                                           ~~\mbox{for some}~ s \in {\cal I}_{\cal F} $ \\

But for $x \in \mathbb{R}^\Lambda$, we also have the equivalent conditions \\

$~~~~~~ x \in {\cal I}_{\cal F} ~~~\Longleftrightarrow~~~ Z ( x ) \in {\cal F} $ \\

and the proof of (3.11) is completed. \\

In view of (3.11), it follows that every ideal in $\mathbb{R}^\Lambda$ is of the form
${\cal I}_{\cal F}$, where ${\cal F}$ is a certain filter on $\Lambda$. Also, every filter on
$\Lambda$ is of the form ${\cal F}_{\cal I}$, where ${\cal I}$ is a certain ideal in
$\mathbb{R}^\Lambda$. \\ \\

{\bf 4. Further large classes of Algebras of Scalars} \\

Related to recent studies of nonlinear PDEs and their generalized solutions, quotient algebras
of scalars of {\it more general form} than  those in (3.4) have been introduced and used in
order, among others, to assign point values to generalized solutions. Such algebras of scalars,
which turn out to be natural extensions of the field $\mathbb{R}$ of real numbers, or of the
field $\mathbb{C}$ of complex numbers, are Non-Archimedean and have zero divisors, thus they
are {\it no longer} division algebras, see Colombeau [1,2], Biagioni, Oberguggenberger,
Grosser, et. al. [1,2], Garetto. \\
For the original and general set up of quotient algebras used in solving large classes of
nonlinear PDEs, see Rosinger [2-15], Rosinger \& Walus [1,2], Mallios \& Rosinger [1,2].
Further details in this regard can be found at the whole corresponding field 46F30 in the
AMS 2000 Subject Classification at www.ams.org/index/msc/46Fxx.html, or Grosser, et. al.
[2, p. 7]. \\

In fact, back in the late 1950s, similarly general quotient algebras were introduced in order
to construct nonstandard extensions of the field $\mathbb{R}$ of real numbers, see Schmieden
\& Laugwitz. This construction of nonstandard reals predates by a few years that of A
Robinson. However, the algebras of Schmieden \& Laugwitz were {\it not} fields, as they had
zero divisors. This is one of the reasons why the construction of the nonstandard field of
reals $^*\mathbb{R}$ by Robinson became the preferential one. \\

For the sake of clarifying the background, we recall in short that the general construction of
the algebras used in the study of nonlinear PDEs proceeds in a rather simple and direct manner,
as follows. Let $\mathbb{K}$ be any algebra, and let us consider \\

(4.1) $~~~~~~ {\cal I} ~\subset~ {\cal A} ~\subseteq~ \mathbb{K}^\Lambda $ \\

with ${\cal A}$ any subalgebra in $K^\Lambda$ and ${\cal I}$ any ideal in ${\cal A}$. Then we
obtain the quotient algebra \\

(4.2) $~~~~~~ A ~=~ {\cal A} / {\cal I} $ \\

Now there is a {\it natural} mapping of $\mathbb{K}$ into these algebras $A$, namely \\

(4.3) $~~~~~~ \mathbb{K} \ni x ~~\longmapsto~~ u ( x ) + {\cal I} \in A ~=~
                                                     {\cal A} / {\cal I} $ \\

where $u ( x ) \in \mathbb{K}^\Lambda$ is the {\it constant} function $u ( x ) : \Lambda
\longrightarrow \mathbb{K}$ for which $u ( x ) ( \lambda ) = x$, for $\lambda \in \Lambda$. Of
course, in order for the mapping (4.3) to be well defined, the quotient algebra (4.2) has to
satisfy the condition \\

(4.4) $~~~~~~ {\cal U}~_\Lambda ~\subseteq~ {\cal A} $ \\

where \\

(4.5) $~~~~~~ {\cal U}~_\Lambda ~=~ \{~ u ( x ) ~~|~~ x \in \mathbb{K} ~\}
$ \\

is the {\it diagonal} of the Cartesian product $\mathbb{K}^\Lambda$, and itself is a
subalgebra of $\mathbb{K}^\Lambda$. \\
Also, it is easy to see that the mapping (4.3) is {\it injective}, that is, as an algebra,
$\mathbb{K}$ is {\it embedded} in $A$, or equivalently, $A$ is an {\it algebra extension} of
$\mathbb{K}$, if and only if \\

(4.6) $~~~~~~ {\cal I}~ \bigcap~ {\cal U}~_\Lambda ~=~ \{~ 0 ~\} $
\\

in other words, the ideal ${\cal I}$ is {\it off diagonal}. \\

Clearly, the conditions (4.4) and (4.6) which characterize the algebra embedding or
extension \\

(4.7) $~~~~~ \mathbb{K} ~\subseteq~ A ~=~ {\cal A} / {\cal I} $ \\

given by the natural mapping (4.3) are easy to satisfy, as seen in the literature mentioned
above. Consequently, there are plenty of such algebra embeddings or extensions. \\

We note that the algebras (3.4) are particular cases of the algebras (4.2), namely,
corresponding to ${\cal A} = \mathbb{K}^\Lambda$. And in such a situation condition (4.4) is
obviously satisfied. \\

In a subsequent paper, we shall perform the computation of Feynman path integrals precisely in
such algebras (4.2), (4.4), (4.6) \\ \\

{\bf 5. On Representation and Interpretation in Mathematics} \\

Before we proceed further, a few additional remarks may be useful. \\

Lack of one's longer time familiarity with certain new mathematical structures can often
constitute a considerable threshold which may easily prevent one from crossing it. \\
This could already be noted during the last century and half, in the case of going beyond the
traditional fields of real or complex scalars, $\mathbb{R}$ and $\mathbb{C}$, respectively,
and working with the quaternions $\mathbb{H}$, let alone, with the octonions $\mathbb{O}$. \\
And one of the reasons for the manifestation of such thresholds is in the lack of a
sufficiently familiar {\it interpretation} of the respective new mathematical structures. \\

Yet, we cannot fail to note that one of the main powers of mathematics rests precisely in the
fact that, so often in mathematical theories, their ability for a rigorous {\it representation}
goes considerably {\it farther} than our ability for a familiar enough {\it interpretation}. \\

And to give a simple, yet particularly relevant example, let us recall that the Peano Axioms
for the natural numbers $\mathbb{N}$ offer a rigorous {\it representation} for all $n \in
\mathbb{N}$. On the other hand, we are - and can only be - familiar enough with the {\it
interpretation} of a few and not too large such $n \in \mathbb{N}$. Indeed, the number

$$ n ~=~ 10^{10^{10}} $$

for instance, only needs six digits for a perfectly rigorous representation. As for its
interpretation, on the other hand, its meaning seems completely to fail us. Certainly, this
$n$ is supposed to be far larger than the total number of elementary particles in the whole of
presently known Universe. Not to mention that we do not have any meaningful, let alone
familiar way of interpretation not relaying on the above representation, and which would alow
us to distinguish between this particular $n$ and, say, its immediate successor $n + 1$. In
fact, in the same manner, we could not possibly distinguish between this $n$ and, say, $n +
1000000$ either. \\

Consequently, in the construction of a mathematical theory one should avoid the situation
where interpretation is unduly restricting representation. Not to mention that, since Cantor's
Set Theory, and of the more recent Category Theory, it becomes more obvious that one of the
main powers of mathematics comes precisely from its ability to produce {\it impressively
general representations}, and do so not constrained by the possible lack of more familiar
interpetations. \\

As for the tendency so far in the study of Feynman path integrals to impose a rather exclusive
use of Functional Analysis, with the respective topological and/or measure structures on
infinite dimensional spaces of functions, structures which by now offer a particularly
familiar world of interpretations, one may recall that Nonstandard Analysis is a good and
recent example of a mathematical theory where interpretation has not been allowed to limit,
let alone, preclude the power of representation. And so far, none of the methods which
constructs the nonstandard reals $^*\mathbb{R}$ makes any use of Functional Analysis. \\

For further details on the possible interplay between representation and interpretation see
Rosinger [8, pp. 391-402], [9, chap. 1, sect. 7]. \\

Needless to say, similar problems of interpretation we may face, when dealing with the
extended algebras of scalars mentioned in section 4. \\
On the other hand, to the extent that such extended algebras of scalars do allow a rigorous
computation of the Feynman path integrals, their use - as means of representation - should not
be set aside so easily. And especially not, in view of the long outstanding lack of any other
alternative rigorous representation. Not to mention the importance of Feynman path integrals
in Quantum Theory. \\ \\

{\bf 6. The need for Non-Archimedean Scalars} \\

The issue of "infinities in physics" has a troubling history. In Quantum Field Theory, it is
attempted to be dealt with by renormalization. \\
Here we shall look into one of the simplest basic mathematical phenomena involved, one that
points to the possible interest in considering algebras of scalars extending the usual field
$\mathbb{C}$ of complex numbers. These algebras, although commutative and associative, may
turn out to have zero divisors, and furthermore, be no longer Archimedean. \\

One of the simplest ways to approach this issue of "infinities" is to see it within the usual
Differential Calculus where it can be treated through the so called Banach generalized limits,
recalled here for convenience in its essence, see Berberian [pp. 117-122], or Paterson, for
the general theory. \\

Let ${\cal B}$ be the set of all usual {\it bounded} sequences of real numbers, that is, the
set of all the bounded functions $x : \mathbb{N} \longrightarrow \mathbb{R}$. With the usual
term-wise operations on such sequences, ${\cal B}$ is a vector space on $\mathbb{R}$. Defining on
${\cal B}$ the norm

$$ | | x | | ~=~ \sup_{n \in \mathbb{N}}~ | x ( n ) |,~~~ x \in {\cal B} $$

we obtain a Banach space structure on ${\cal B}$. \\
Let ${\cal C}onv$ be the subset of all usual {\it convergent} sequences in ${\cal B}$. This
means that there exists a {\it unique} linear functional \\

(6.1) $~~~~~~ L : {\cal C}onv ~\longrightarrow~ \mathbb{R} $ \\

such that for every convergent sequence $x \in {\cal C}onv$, we have \\

(6.2) $~~~~~~ L ( x ) ~=~ \lim_{n \to \infty}~ x ( n ) $ \\

In other words, $L$ is the usual limit of convergent sequences. As is well known, this linear
functional $L$ has the following properties \\

(6.3) $~~~~~~ \inf_{n \in \mathbb{N}}~ x ( n ) ~\leq~ L ( x ) ~\leq~
                                  \sup_{n \in \mathbb{N}}~ x ( n ),~~~ x \in {\cal C}onv $ \\

(6.4) $~~~~~~ L ( x^{+} ) ~=~ L ( x ),~~~ x \in {\cal C}onv $ \\

where $x^{+}$ is the sequence obtained from $x$ by one {\it shift} to the right, namely,
$x^{+} ( n ) = x ( n + 1 )$, for $n \in \mathbb{N}$. \\

The nontrivial fact is that the linear functional $L$ which is defined on ${\cal C}onv$ by
(6.2), can be {\it extended} to a linear functional $L^{\#}$ on the whole of ${\cal B}$, with
the {\it preservation} of the above properties (6.3) and (6.4), Berberian. \\
However, this extension is {\it not} unique. One of the reasons for this lack of uniqueness is
in the celebrated Hahn-Banach theorem in Functional Analysis about the prolongation of linear
functionals, a theorem used essentially in order to obtain the extension $L^{\#}$ of $L$. \\

The above can be rewritten equivalently in terms of usual {\it series} of real numbers. Indeed,
given any such series \\

(6.5) $~~~~~~ \Sigma_{n \in \mathbb{N}}~ s ( n ) $ \\

where $s : \mathbb{N} \longrightarrow \mathbb{R}$, we can associate with it uniquely the
sequence $x : \mathbb{N} \longrightarrow \mathbb{R}$, defined by \\

(6.6) $~~~~~~ x ( n ) ~=~ s ( 0 ) + s ( 1 ) + s( 2 ) + ~.~.~.~ + s ( n ),~~~
                                                             n \in \mathbb{N} $ \\

Then, as is well known, we have by definition \\

(6.7) $~~~~~~ \Sigma_{n \in \mathbb{N}}~ s ( n ) ~=~ S
                                  ~~~\Longleftrightarrow~~~ L ( x ) ~=~ S $ \\

in other words, the series (6.5) is convergent and converges to the real number $S \in
\mathbb{R}$ which is its sum, if and only if the sequence (6.6) is convergent and converges to
the same real number. It follows that \\

(6.8) $~~~~~~ \begin{array}{l}
                  \Sigma_{n \in \mathbb{N}}~ s ( n )~~ \mbox{is divergent}
                  ~~~\Longleftrightarrow~~~ x \in \mathbb{R}^\mathbb{N} \setminus {\cal C}onv
                  ~~~\Longleftrightarrow~~~ \\ \\
                  ~~~\Longleftrightarrow~~~ x \in {\cal B} \setminus {\cal C}onv
                          ~~\mbox{or}~~ x \in \mathbb{R}^\mathbb{N} \setminus {\cal B}
               \end{array} $ \\

Now, in terms of series such as $\Sigma_{n \in \mathbb{N}}~ s ( n )$ in (6.5), we encounter
"infinities" precisely when we have $x \in \mathbb{R}^\mathbb{N} \setminus {\cal B}$, for the
corresponding sequences in (6.6). A relevant study in this respect can be found in Hardy,
which happens to be the last more important mathematical contribution, published posthumously,
of that noted English mathematician. \\

Consequently, those series $\Sigma_{n \in \mathbb{N}}~ s ( n )$ which correspond to
"infinities" can no longer be treated either by the linear functional $L$ on ${\cal C}onv$, or
by its extension $L^{\#}$ to ${\cal B}$, since the corresponding sequences $x$ are {\it no
longer} bounded. \\

Here it is relevant to note that the linear functional $L^{\#}$ {\it cannot} further be
extended much beyond ${\cal B}$, as long as we want to preserve the above property (6.4). \\
Indeed, let us assume that there exists a linear functional \\

(6.9) $~~~~~~ L^{\#\#} : \mathbb{R}^\mathbb{N} ~\longrightarrow~ \mathbb{R} $ \\

which extends $L^{\#}$ and which satisfies (6.4). We take then the {\it unbounded} sequence
$\nu \in \mathbb{R}^\mathbb{N}$ defined by $\nu ( n ) = n$, for $n \in \mathbb{N}$. In this
case, it is obvious that

$$ \nu^{+} ~-~ \nu ~=~ ( 1, 1, 1, ~.~.~.~ ) \in {\cal C}onv $$

hence we obtain

$$ L^{\#\#} ( \nu^{+} ) ~-~ L^{\#\#} ( \nu ) ~=~ L^{\#\#} ( \nu^{+} ~-~ \nu ) ~=~ $$
$$ ~~~~~~~=~ L^{\#\#} ( ( 1, 1, 1, ~.~.~.~ ) ) ~=~ L ( ( 1, 1, 1, ~.~.~.~ ) ) ~=~ 1 $$

However, in view of (6.4), we have $L^{\#\#} ( \nu^{+} ) = L^{\#\#} ( \nu )$, thus it follows
that $0 = 1$, and we obtained a contradiction. \\

About the above proof it is important to note that one does not in fact need the strong
assumption in (6.9), namely that $L^{\#\#}$ be defined on the whole of
$\mathbb{R}^\mathbb{N}$. Indeed, the above contradiction follows even if we only assume that
$L^{\#\#}$ is extended beyond ${\cal B}$ to some larger vector subspace ${\cal E}$ of
$\mathbb{R}^\mathbb{N}$, provided that we have $\nu \in {\cal E}$. \\

Also it is important to note that the {\it unbounded} sequence $\nu$ which prevents the
further extension of the linear functional $L^{\#}$ does {\it not} grow fast to infinity. \\

Let us now summarize the above. \\

Series which are connected to "infinities" correspond to {\it unbounded} sequences. \\
The Banach generalized limits {\it cannot} be extended even to moderately growing unbounded
sequences. \\

And here we come to the {\it likelihood} of having to deal with non-Archimedean algebras of
scalars. Indeed, in order to be able to deal with "infinities", the quotient algebras (3.4),
or more generally (4.2), will often contain scalars of the form

$$ a ~=~ x ~+~ {\cal I},~~~ x \in \mathbb{K}^\Lambda, ~\mbox{respectively},~ x \in {\cal A} $$

where $x$ are unbounded. As a result, such algebras tend not to be Archimedean. \\
A good illustration of this phenomenon happens in the particular case of the algebras (3.4)
which give the field $^*\mathbb{R}$ of nonstandard reals. \\ \\

{\bf 7. Conclusions} \\

Since the introduction of Feynman path integrals in the 1940s, since the introduction of
various methods of "re-normalization", or since the need to deal with various "infinities in
physics", the ongoing unsuccessful efforts to find for them a sufficiently rigorous and
systematic mathematical approach have been characterized by {\it two restrictive assumptions},
both of them accepted as self-evident, and never questioned. Namely :

\begin{itemize}

\item the Feynman path integrals, the results of re-normalization, as well as of other
infinities in Physics must by all means have values given by real or complex numbers,

\item Functional Analysis, and specifically, topology and/or measure theory on suitable
infinite dimensional spaces of functions is the only way to make the Feynman path integrals
mathematically rigorous. As for re-normalization or other infinities in Physics, any mix of
existing mathematical gimmicks, among those well known to physicists, may do.

\end{itemize}

In this paper both above restrictive assumptions are set aside. Instead, it is shown that one
can easily construct a large class of algebras of scalars which extend the complex numbers. \\
These algebras of scalars have two features : they are no longer Archimedean, and they have
zero divisors. However, they are both commutative and associative. \\
The loss of the Archimedean property turns out to present certain advantages, since it helps
in avoiding the so called "infinities in physics". \\
As for the presence of zero divisors, we have for long been accustomed to this phenomenon,
ever since we have been dealing with algebras of usual matrices of real or complex
numbers. \\

Here it should be mentioned that a foundational type questioning of certain assumptions which
are usually taken for granted in physics, and which may have connections with the above, has
occurred in the literature, see Litvinov et.al, or Rosinger [18]. \\

In a subsequent paper we shall present the way the Feynman path integrals can at last be given
{\it rigorously} computed values in such algebras of scalars. \\
In a further paper, we shall redefine the Hilbert spaces in terms of such algebras of scalars.
Hilbert spaces, as known, are fundamental in Quantum Mechanics, yet back in 1935, John von
Neumann, in a letter to George David Birkhoff, wrote : "I WOULD LIKE TO MAKE A CONFESSION
WHICH MAY SEEM IMMORAL : I DO NOT BELIEVE IN HILBERT SPACE ANYMORE." Some of the reasons for
that confession may have been simple but basic failures, such as the fact that neither the
position, nor the momentum operators have eigenvectors in Hilbert spaces.  \\ \\

\end{document}